\ifx\shlhetal\undefinedcontrolsequence\let\shlhetal\relax\fi

\documentclass[11pt]{amsart}

\usepackage{amsmath}
\usepackage{amssymb}

\newtheorem{theorem}{Theorem}[section]
\newtheorem{claim}[theorem]{Claim}

\theoremstyle{definition}
\newtheorem{definition}[theorem]{Definition}

\newtheorem{question}[theorem]{Question}

\theoremstyle{remark}

\newcount\skewfactor
\def\mathunderaccent#1#2 {\let\theaccent#1\skewfactor#2
\mathpalette\putaccentunder}
\def\putaccentunder#1#2{\oalign{$#1#2$\crcr\hidewidth
\vbox to.2ex{\hbox{$#1\skew\skewfactor\theaccent{}$}\vss}\hidewidth}}
\def\name{\mathunderaccent\tilde-3 }

\def\smallbox#1{\leavevmode\thinspace\hbox{\vrule\vtop{\vbox
   {\hrule\kern1pt\hbox{\vphantom{\tt/}\thinspace{\tt#1}\thinspace}}
   \kern1pt\hrule}\vrule}\thinspace}

\newcommand{\cf}{{\rm cf}}

\usepackage[all]{xy}

\def\qedref#1{$\qed_{\reforiginal{#1}}$}

\title{On the spectrum of characters of ultrafilters}
\author{Shimon Garti}
\address{Institute of Mathematics,
 The Hebrew University of Jerusalem,
 Jerusalem 91904, Israel}
\email{shimon.garty@mail.huji.ac.il}

\author{Menachem Magidor}
\address{Institute of Mathematics,
 The Hebrew University of Jerusalem,
 Jerusalem 91904, Israel}
\email{menachem@math.huji.ac.il}

\author{Saharon Shelah}
\address{Institute of Mathematics,
 The Hebrew University of Jerusalem,
 Jerusalem 91904, Israel,
 and Department of Mathematics
 Rutgers University
 New Brunswick, NJ 08854, USA}
\email{shelah@math.huji.ac.il}
\urladdr{http://www.math.rutgers.edu/\char`\~shelah}

\thanks{Second author's research was partially supported by the ISF grant no. 817/11. Third author's research was partially supported by the United States-Israel Binational Science Foundation, and ISF grant no. 1053/11. This is publication 1040 of the third author.}

\begin{document}
\let\labeloriginal\label
\let\reforiginal\ref

\begin{abstract}
We show that the character spectrum ${\rm Sp}_\chi(\lambda)$ (for a singular cardinal $\lambda$ of countable cofinality) may include any prescribed set of regular cardinals between $\lambda$ and $2^\lambda$. \newline 
Nous prouvons que ${\rm Sp}_\chi(\lambda)$ (par un cardinal singulier $\lambda$ avec cofinalit\`e nombrable) peut comporter tout l'ensemble prescrit de cardinaux reguliers entre $\lambda$ et $2^\lambda$.
\end{abstract}

\maketitle

\newpage

\section{introduction}

The purpose of this paper is to analize the spectrum of characters of uniform ultrafilters. Recall that an ultrafilter $U$ on $\lambda$ is called uniform if the size of every member of $U$ is $\lambda$. We need some basic definitions for describing the contents of this paper:

\begin{definition}
\label{ccharacters} Characters of ultrafilters. \newline 
Let $\lambda$ be an infinite cardinal, and $U$ an ultrafilter on $\lambda$.
\begin{enumerate}
\item [$(\aleph)$] A base $\mathcal{A}$ for $U$ is a subfamily of $U$ such that for every $B\in U$ there is some $A\in\mathcal{A}$ with the property $A\subseteq^*B$.
\item [$(\beth)$] The character of $U$, denoted by ${\rm Ch}(U)$, is the minimal cardinality of a base for $U$.
\item [$(\gimel)$] ${\rm Sp}_\chi(\lambda)$ is the set of all $\mu$-s so that $\mu={\rm Ch}(U)$ for some uniform ultrafilter $U$ on $\lambda$.
\item [$(\daleth)$] The ultrafilter number $\mathfrak{u}_\lambda$ is the minimal value of ${\rm Ch}(U)$ for some uniform ultrafilter on $\lambda$.
\end{enumerate}
\end{definition}

It is known that $\mathfrak{u}_\lambda>\lambda$ for every infinite cardinal $\lambda$ (see, e.g., \cite{MR2992547}, Claim 1.2). It follows that if $2^\lambda=\lambda^+$ then ${\rm Ch}(U)=\lambda^+$ for every uniform ultrafilter $U$. A natural question is whether $\mathfrak{u}_\lambda=\lambda^+$ is consistent with large values of $2^\lambda$. A positive answer is supplied in \cite{MR2992547}. It is proved there that $\mathfrak{u}_\lambda=\lambda^+$ is consistent with arbitrarily large value of $2^\lambda$, for some singular cardinal $\lambda$.

The proof requires a singular cardinal, limit of measurables. Let us try to describe the philosophy of the proof. We begin with a supercompact cardinal $\lambda$, making it singular via Prikry forcing or Magidor forcing. As mentioned above, an easy way to produce a uniform ultrafilter with a small generating base is to invoke the continuum hypothesis. But how can we enlarge $2^\lambda$ while keeping this special ultrafilter?

Well, we still employ the continuum hypothesis, but this time we apply it to a sequence of measurable cardinals below $\lambda$. We choose such a sequence, $\langle\lambda_i:i<\cf(\lambda)\rangle$, so that $2^{\lambda_i}=\lambda_i^+$ for every $i<\cf(\lambda)$. We choose a normal ultrafilter $U_i$ on each $\lambda_i$, and by virtue of the continuum hypothesis on $\lambda_i$ we have a small base for $U_i$. Now we combine these bases in such a way that enables us to elicit a base for an ultrafilter $U$ on $\lambda$ of size $\lambda^+$.

The main point here is that the base of $U$ is created from the small bases of the $U_i$-s, hence there is no need to diminish $2^\lambda$. The assumption $2^\lambda=\lambda^+$ can be replaced by the continuum hypothesis on each $\lambda_i$. More precisely, we need $2^{\lambda_i}=\lambda_i^+$ for every $i<\cf(\lambda)=\kappa$ and some kind of approachability of $\lambda^+$. It suffices that ${\rm tcf}(\prod\limits_{i<\kappa}\lambda_i,<_{J_\kappa^{\rm bd}})={\rm tcf}(\prod\limits_{i<\kappa}\lambda_i^+,<_{J_\kappa^{\rm bd}})=\lambda^+$. Under these assumptions, the desired $U$ can be created even if $2^\lambda$ is large.

Our problem invades now to the realm of \emph{pcf} theory. Is it possible to have a singular cardinal $\lambda$, limit of measurables, such that the continuum hypothesis holds on these measurables, the true cofinalities of the members of the sequence (as well as their successors) are $\lambda^+$, yet $2^\lambda$ is large? The answer is positive, as shown in \cite{MR2992547}.

The forcing machinery employed in \cite{MR2992547} is taken from \cite{MR2987137}. The main component of this forcing is an iteration of `product dominating real', and the last component can be Prikry forcing or Magidor forcing. It follows that the cofinality of $\lambda$ can be any regular $\kappa$ below it. In \cite{MR2992547} we are seeking for $\lambda^+$ as the value of ${\rm Ch}(U)$. But \cite{MR2987137} provides us with the ability to force ${\rm Ch}(U)=\mu$ for every regular cardinal $\mu$ above $\lambda$. This is done here, and it gives any single value as a possible member in ${\rm Sp}_\chi(\lambda)$. Our way to controll the true cofinalities is simply by determining the length of the iteration in the product dominating real forcing.

Now we would like to add two distinct members $\mu_0$ and $\mu_1$ (or more generally, any prescribed set of regular cardinals), to ${\rm Sp}_\chi(\lambda)$. The na\"ive approach is to designate two distinct sequences of measurables below $\lambda$. One sequence would catch $\mu_0$, and the other would try to catch $\mu_1$. But here we are confronted with an obstacle.

If we just concatenate two iterations, each of the desired length (aiming to seize both $\mu_0$ and $\mu_1$), we fail. The second iteration ruins the achievements of the first one. Hence we must find a different way to controll distinct sequences of measurables and their true cofinalities.

For this, we employ in this paper the Extender-based Prikry forcing. We begin with a strong cardinal $\lambda$, and a large set of measurables above it. The basic step blows up $2^\lambda$ to any desired point, and creates many distinct sequences of measurables with different values of true cofinalities. For every measurable $\mu$ above $\lambda$ we can find a sequence of measurables below $\lambda$ so that ${\rm tcf}(\prod\limits_{i<\kappa}\lambda_i,<_{J_\kappa^{\rm bd}})=\mu$ and ${\rm tcf} (\prod\limits_{i<\kappa}\lambda_i^+,<_{J_\kappa^{\rm bd}})=\mu^+$.
It should be emphasized that we preserve the full GCH below $\lambda$. Consequently, we must confine ourselves to a singular $\lambda$ with countable cofinality, since $2^\lambda$ becomes large (and due to Silver's theorem).

This basic step presents every successor of measurable above $\lambda$ as an element in ${\rm Sp}_\chi(\lambda)$. The next step is to realize each regular cardinal $\chi$ (above $\lambda$) as a member of ${\rm Sp}_\chi(\lambda)$. For this end, we collapse each successor of the measurable cardinal $\mu_i$ to the regular cardinal $\chi_i$ while keeping the true cofinalities. At the end, every set of regular cardinals in the interval $(\lambda,2^\lambda]$ can be realized as a set of values for ${\rm Sp}_\chi(\lambda)$.

Let us indicate that almost full caharcterization of ${\rm Sp}_\chi(\lambda)$ appears in \cite{Sh915} for the case of $\lambda=\aleph_0$. Such a characterization requires two orthogonal methods. The first one is a method for inserting a member into ${\rm Sp}_\chi(\lambda)$, and the second is a method for eliminating the memberhood of a cardinal. In our context, we supply a general way for inserting, but we (still) do not know how to eliminate.

At last, one has to admit the existence of some large cardinals in the ground model. But inasmuch as $\aleph_0$ is a large cardinal, it seems that there is no philosophical reason to deny the existence of other large cardinals.

Our notation is standard. We refer to \cite{MR1318912} and \cite{MR2768693} for \emph{pcf} theory, \cite{MR2768685} for general background on cardinal invariants, and \cite{MR2768695} for the subject of Prikry-type forcings. Our notation is coherent, in general, with these monographs. In particular, we adopt the Jerusalem notation in forcing, i.e., $p\leq q$ means that $q$ is a stronger condition than $p$.

The relation $A\subseteq^* B$ reads $A$ is almost included in $B$. If $A$ and $B$ are subsets of $\lambda$ then $A\subseteq^* B$ means that $|A\setminus B|<\lambda$. We use $\theta,\kappa,\lambda,\mu,\chi$ for infinite cardinals, and $\alpha,\beta,\gamma,\delta,\varepsilon,\zeta$ as well as $i,j$ for ordinals.
For a sequence of cardinals $\bar{\lambda}=\langle\lambda_i:i<\kappa\rangle$ we denote $\bigcup\{\lambda_j:j<i\}$ by $\lambda_{<i}$. If $\bar{\lambda}$ is the sequence $\langle\lambda_i:i<\kappa\rangle$ then $\bar{\lambda}^+$ denotes the sequence of the successors, i.e., $\langle\lambda_i^+:i<\kappa\rangle$. For a regular cardinal $\kappa$, ${J_\kappa^{\rm bd}}$ is the ideal of bounded subsets of $\kappa$.

\newpage

\section{The spectrum of characters}

Let us begin with the product dominating real forcing $\mathbb{Q}_{\bar{\theta}}$. We start with a Laver-indestructible supercompact cardinal $\lambda$, assuming that the GCH holds above $\lambda$. We choose an increasing sequence of regular cardinals $\bar{\theta}$, which steps up fast enough. The following is the basic component of the forcing notion:

\begin{definition}
\label{dominating}
The $\bar{\theta}$-dominating forcing. \newline 
Let $\lambda$ be a supercompact cardinal.
Suppose $\bar{\theta} = \langle \theta_\alpha : \alpha < \lambda \rangle$ is an increasing sequence of regular cardinals so that $2^{|\alpha|+\aleph_0} < \theta_\alpha < \lambda$ for every $\alpha < \lambda$. 
\begin{enumerate}
\item [$(\aleph)$] $p \in \mathbb{Q}_{\bar{\theta}}$ iff:
\begin{enumerate}
\item $p = (\eta, f) = (\eta^p, f^p)$,
\item $\ell g(\eta) < \lambda$,
\item $\eta \in \prod \{ \theta_\zeta : \zeta < \ell g(\eta) \}$,
\item $f \in \prod \{ \theta_\zeta : \zeta < \lambda \}$,
\item $\eta \triangleleft f$ (i.e., $\eta(\zeta)=f(\zeta)$ for every $\zeta<\ell g(\eta)$).
\end{enumerate}
\item [$(\beth)$] $p \leq_{\mathbb{Q}_{\bar{\theta}}} q$ iff ($p,q \in \mathbb{Q}_{\bar{\theta}}$ and) 
\begin{enumerate}
\item $\eta^p \trianglelefteq \eta^q$,
\item $f^p(\varepsilon) \leq f^q(\varepsilon)$, for every $\varepsilon < \lambda$.
\end{enumerate}
\end{enumerate}
\end{definition}

We iterate this forcing, and then we compose Prikry forcing (or Magidor forcing) to make $\lambda$ a singular cardinal of cofinality $\kappa$. If we iterate the product dominating real forcing along some ordinal $\delta\geq\lambda^{++}$ such that $|\delta|=|\delta|^\lambda$ then $2^\lambda$ becomes $|\delta|$ and ${\rm tcf}(\prod\limits_{i<\kappa}\theta_i,<_{J_\kappa^{\rm bd}})={\rm tcf}(\prod\limits_{i<\kappa}\theta_i^+,<_{J_\kappa^{\rm bd}})=\cf(\delta)$. The following theorem (from \cite{MR2987137}) summarizes the pertinent properties:

\begin{theorem}
\label{ffforcing} Product dominating $\lambda$-reals. \newline 
Assume there is a supercompact cardinal in the ground model. \newline 
Then one can force the existence of a singular cardinal $\lambda>\cf(\lambda)=\kappa$, a limit of measurables $\bar{\lambda}=\langle\lambda_i:i<\kappa\rangle$, such that $2^{\lambda_i}=\lambda_i^+$ for every $i<\kappa, 2^\lambda\geq |\delta|$ and the true cofinality of both products $\prod\limits_{i<\kappa}\lambda_i /J_\kappa^{\rm bd}$ and $\prod\limits_{i<\kappa}\lambda^+_i /J_\kappa^{\rm bd}$ is $\cf(\delta)$ for some prescribed $\delta\geq\lambda^{++}$.
\end{theorem}

\hfill \qedref{ffforcing}

We employ the above theorem in order to prove that each regular $\mu$ above $\lambda$ may be a member of ${\rm Sp}_\chi(\lambda)$:

\begin{claim}
\label{oonestep} Forcing one member. \newline 
Suppose there is a supercompact cardinal $\lambda$ in the ground model, and $\mu$ is a regular cardinal above $\lambda$. \newline 
Then one can force $\mu\in{\rm Sp}_\chi(\lambda)$, upon making $\lambda$ a singular cardinal.
\end{claim}

\par\noindent\emph{Proof}.\newline 
Due to Theorem \ref{ffforcing}, we may assume that $\lambda>\cf(\lambda)=\kappa$ is the limit of a sequence of measurable cardinals $\bar{\lambda}= \langle\lambda_i:i<\kappa\rangle$, and the continuum hypothesis holds for every $\lambda_i$ (upon using the product dominating real forcing). We choose a normal ultrafilter $U_i$ on $\lambda_i$ for every $i<\kappa$, and a fixed uniform ultrafilter $E$ on $\kappa$.

Let $\langle A_{i,\alpha}:\alpha<\lambda_i^+\rangle$ be a $\subseteq^*$-decreasing sequence of sets, which serves as a base for $U_i$ for every $i<\kappa$. Here we use the normality of $U_i$. We choose any base of size $\lambda^+$ to $U_i$ and enumerate its members by $\langle B_\alpha:\alpha<\lambda^+\rangle$. Now we replace each $B_\alpha$ (apart from the first one, $B_0$) by the diagonal intersection of all the previous sets.
By Theorem \ref{ffforcing} we can arrange $2^\lambda\geq\mu$ and ${\rm tcf}(\prod\limits_{i<\kappa}\lambda_i,<_E)= {\rm tcf}(\prod\limits_{i<\kappa}\lambda_i^+,<_E)=\mu$. Let $\bar{f}= \langle f_\alpha:\alpha<\mu\rangle$ be cofinal in $(\prod\limits_{i<\kappa}\lambda_i,<_E)$ and $\bar{g}= \langle g_\beta:\beta<\mu\rangle$ be cofinal in $(\prod\limits_{i<\kappa}\lambda_i^+,<_E)$.

For each $\alpha<\mu,\beta<\mu$ and $Y\in E$ we define the following set:

$$
B_{\alpha,\beta,Y}= \{\zeta<\lambda:\exists i\in Y, \zeta\in[\lambda_{<i}, \lambda_i) \wedge \zeta\in A_{i,g_\beta(i)}\setminus f_\alpha(i)\}.
$$

We collect these sets into $\mathcal{B}=\{B_{\alpha,\beta,Y}:\alpha,\beta<\mu, Y\in E\}$. The cardinality of $\mathcal{B}$ is $\mu$, and it serves as a base for a uniform ultrafilter $U$ on $\lambda$ (see the proof of the main theorem in \cite{MR2992547}), hence ${\rm Ch}(U)\leq\mu$. Let us show that equality holds.

Assume towards contradiction that $\mathcal{C}$ is another base for $U$, and $|\mathcal{C}|=\mu'<\mu$. Without loss of generality $\mathcal{C}\subseteq \mathcal{B}$ (as one can replace each member of $\mathcal{C}$ by an element from $\mathcal{B}$, since $\mathcal{B}$ is also a base for $U$). So $\mathcal{C}$ is a collection of sets of the form $B_{\alpha,\beta,Y}$.

Let $\beta_*$ be an ordinal which is larger than every $\beta$ mentioned in any $B_{\alpha,\beta,Y}\in\mathcal{C}$ but smaller than $\mu$. The collection $\mathcal{A}_{i,\beta_*}= \{A_{i,\gamma}:\gamma<g_{\beta_*}(i)\}$ is not a base of $U_i$ (for any $i<\kappa$). Choose $X_i\subseteq\lambda_i$ so that neither $X_i$ nor $\lambda_i\setminus X_i$ is generated by $\mathcal{A}_{i,\beta_*}$ (i.e., for every $A\in\mathcal{A}_{i,\beta_*}$ we have $|A\setminus X_i|=|A\cap X_i|=\lambda_i$).

Set $X=\bigcup\{X_i:i<\kappa\}$, and assume towards contradiction that $(B_{\alpha,\beta,Y}\subseteq^* X)\vee(B_{\alpha,\beta,Y}\subseteq^* \lambda\setminus X)$ for some $B_{\alpha,\beta,Y}\in\mathcal{C}$. Without loss of generality $B_{\alpha,\beta,Y}\subseteq^* X$ (the opposite option is just the same). We may assume that $B_{\alpha,\beta,Y}$ is $\bigcup\{ A_{i,g_\beta(i)}: i<\kappa\}$ (i.e., we assume that $Y=\kappa$ and $f_\alpha(i)=0$ for every $i<\kappa$).

Let $A'_i$ be $A_{i,g_\beta(i)}\setminus X_i$ for every $i<\kappa$. The cardinality of $A'_i$ is $\lambda_i$, as $A_{i,g_\beta(i)}\in\mathcal{A}_{i,\beta_*}$ and $g_{\beta}(i)<g_{\beta_*}(i)$ (without loss of generality, for every $i<\kappa$). It follows that $|\bigcup\{A'_i:i<\kappa\}|=\lambda$. Since $\bigcup\{A'_i:i<\kappa\}\subseteq \lambda\setminus X$ we have $|B_{\alpha,\beta,Y}\setminus X|=\lambda$, a contradiction.

\hfill \qedref{oonestep}

The main forcing notion to be used in this paper is the extender-based Prikry forcing. We follow the notation of \cite{MR2768695} and we refer to the theorems proved there about the extender-based Prikry forcing. Some preliminary definitions and facts are in order.

Let $\kappa$ be a $\lambda$-strong cardinal (i.e., there is an elementary embedding $\jmath:{\rm \bf V}\rightarrow M, \kappa={\rm crit}(\jmath), {\rm V}_\lambda\subseteq M$ and $\jmath(\kappa)>\lambda$). For every $\alpha<\lambda$ we define an ultrafilter $U_\alpha$ on $\kappa$ as follows:

$$
A\in U_\alpha \Leftrightarrow \alpha\in\jmath(A).
$$

The idea is to generate many Prikry sequences in $\kappa$, whence $U_\alpha$ is related to the $\alpha$-th sequence for every $\alpha<\lambda$. Consequently, $2^\kappa\geq\lambda$ in the extension model. Moreover, for every regular cardinal $\alpha<\lambda$ the sequences for all the $\beta$-s below $\alpha$ would be cofinal in the appropriate product (hence the pertinent requirement on the cofinalities is fulfilled).

We begin with preparing a nice system of ultrafilters and embeddings, in order to define the forcing with them. This system is denoted by $E$, and called an extender. For every $\alpha<\lambda$ we have the following commutative diagram:

$$
\xymatrix{
{\rm \bf V} \ar[ddr]_{\imath_\alpha} \ar[r]^\jmath & M \\ \\ 
& {{\rm \bf V}^\kappa/U_\alpha\cong N_\alpha} \ar[uu]_{k_\alpha} }
$$

We use the above fixed elementary embedding $\jmath$, and the canonical embedding $\imath_\alpha$ of ${\rm \bf V}$ into ${\rm \bf V}^\kappa /U_\alpha$. As each $U_\alpha$ is $\kappa$-complete, ${\rm \bf V}^\kappa /U_\alpha$ is well-founded, so one can collapse it onto a transitive model $N_\alpha$ (and we do not distinguish ${\rm \bf V}^\kappa /U_\alpha$ from $N_\alpha$). The mapping $k_\alpha$ is defined by $k_\alpha([f])=\jmath(f) (\alpha)$ (for every equivalence class $[f]\in{\rm \bf V}^\kappa /U_\alpha$). It is routine to check that the diagram is commutative.

Let us define a partial order on the ordinals of $\lambda$. We say that $\beta\leq_E\alpha$ iff ($\beta\leq\alpha$ and) there exists a function $f:\kappa\rightarrow\kappa$ so that $\jmath(f)(\alpha)=\beta$. Intuitively, it means that $\beta$ belongs to the range of $k_\alpha$. Notice that $\leq_E$ is transitive (by compsition of the pertinent functions). It is well known that $\leq_E$ is $\kappa^{++}$-directed.

For each pair of ordinals $(\beta,\alpha)$ such that $\beta\leq_E\alpha$ we choose a projection $\pi_{\alpha\beta}:\kappa \rightarrow\kappa$ which satisfies $\jmath(\pi_{\alpha\beta})(\alpha)=\beta$. Notice that $U_\beta \leq_{\rm RK}U_\alpha$ (as demonstrated by $\pi_{\alpha\beta}$), i.e., $B\in U_\beta \Leftrightarrow \pi_{\alpha\beta}^{-1}(B)\in U_\alpha$. The collection of ultrafilters and projections is our extender $E$, and it enables us to define our forcing notion $\mathbb{P}$.

For adding many Prikry sequences we have at each condition $p\in\mathbb{P}$ a set of ordinals $g=g^p$ which is called the support of $p$ and denoted by ${\rm supp}(p)$. For each $\alpha\in g$ we asign a Prikry sequence, and $g$ contains a $\leq_E$-maximal element for which we hold a $U$-tree. Recall that for an ultrafilter $U$ on $\kappa$ we call $T\subseteq[\kappa]^{<\omega}$ a $U$-tree if $\eta\in T\Rightarrow {\rm Suc}_T(\eta)=\{\alpha<\kappa:\eta^{\frown}(\alpha)\in T\}\in U$.

The last concept that we need is the idea of a permitted ordinal. For every $\alpha\in[\kappa,\lambda)$ and $\nu<\kappa$ let $\nu^\circ$ be the ordinal $\pi_{\alpha\kappa}(\nu)$. The sequence $\langle\nu_0,\ldots,\nu_{n-1}\rangle$ is $\circ$-increasing if $\nu_0^\circ<\ldots<\nu_{n-1}^\circ$. An ordinal $\nu<\kappa$ is permitted for $\langle\nu_0,\ldots,\nu_{n-1}\rangle$ if $\nu^\circ>{\rm max}\{\nu_j^\circ:j<n\}$.

\begin{definition}
\label{ggggitik} The extender-based Prikry forcing. \newline 
Let $\kappa$ be a strong cardinal and $\lambda=\cf(\lambda)>\kappa^+$. We define the forcing notion $\mathbb{P}$ as follows: \newline 
A condition $p\in\mathbb{P}$ is the set $\{\langle\gamma,p^\gamma\rangle: \gamma\in g\setminus\{{\rm max}(g)\}\}\cup \{\langle{\rm max}(g),p^{{\rm max}(g)},T\rangle\}$ such that:
\begin{enumerate}
\item $g\subseteq\lambda, |g|\leq\kappa, {\rm max}(g)$ is a maximal element with respect to $\leq_E$, denoted by ${\rm mc}(p)$.
\item $p^\gamma$ is a finite $\circ$-increasing sequence of ordinals below $\kappa$, for every $\gamma\in g$.
\item $T$ is a $U_{{\rm mc}(p)}$-tree whose trunk is $p^{\rm mc}$ and its members are finite $\circ$-increasing sequences such that $p^{\rm mc}\leq_T \eta_0\leq_T \eta_1 \Rightarrow {\rm Suc}_T(\eta_1)\subseteq {\rm Suc}_T(\eta_0)$.
\item $\pi_{{\rm mc}(p)\gamma}({\rm max}(p^{\rm mc}))$ is not permitted for $p^\gamma$ whenever $\gamma\in g$.
\item $|\{\gamma\in g: \nu$ is not permitted for $p^\gamma\}|\leq \nu^\circ$ for every $\nu\in {\rm Suc}_T (p^{\rm mc})$.
\item $\pi_{{\rm mc}(p)0}$ projects $p^{\rm mc}$ onto $p^0$.
\end{enumerate}
\end{definition}

Let us add the definition of the order (as well as the pure order):

\begin{definition}
\label{oooorder} The forcing order and the pure order. \newline 
Assume $p,q\in\mathbb{P}$. \newline 
We say that $p\leq q$ iff:
\begin{enumerate}
\item ${\rm supp}(p)\subseteq {\rm supp}(q)$.
\item $\gamma\in{\rm supp}(p)\Rightarrow p^\gamma\trianglelefteq q^\gamma$.
\item $q^{{\rm mc}(p)}\in T^p$.
\item If $\gamma\in{\rm supp}(p)$ and $i\in {\rm dom}(q^{{\rm mc}(p)})$ is the largest ordinal such that $q^{{\rm mc}(p)}(i)$ is not permitted for $p^\gamma$ then $q^\gamma\setminus p^\gamma=\pi_{{\rm mc}(p)\gamma}''((q^{{\rm mc}(p)}\setminus p^{{\rm mc}(p)}) \upharpoonright lg(q^{\rm mc}\setminus (i+1)))$.
\item $\pi_{{\rm mc}(q){\rm mc}(p)}$ projects $T^q_{q^{\rm mc}}$ into $T^p_{p^{\rm mc}}$.
\item For every $\gamma\in{\rm supp}(p)$ and $\nu\in{\rm Suc}_{T^q}(q^{\rm mc})$, if $\nu$ is permitted for $q^\gamma$ then $\pi_{{\rm mc}(q)\gamma} = \pi_{{\rm mc}(p)\gamma}(\pi_{{\rm mc}(q){\rm mc}(p)}(\nu))$.
\end{enumerate}
We say that $p\leq^* q$ iff $p\leq q$ and $p^\gamma=q^\gamma$ for every $\gamma\in{\rm supp}(p)$.
\end{definition}

Now we can prove the following theorem:

\begin{theorem}
\label{eeextender} Conrtolling true cofinalities. \newline 
Suppose there is a strong cardinal $\lambda$ and $\langle\mu_i:i\leq i(*)\rangle$ is an increasing sequence of measurables above $\lambda$ in the ground model. Then there is a forcing notion $\mathbb{P}$ so that the following hold in the forcing extension ${\rm \bf V}^{\mathbb{P}}$:
\begin{enumerate}
\item [$(a)$] no cardinal is collapsed by $\mathbb{P}$.
\item [$(b)$] $\lambda>\cf(\lambda)=\aleph_0$.
\item [$(c)$] $\theta<\lambda\Rightarrow 2^\theta=\theta^+$.
\item [$(d)$] $2^\lambda\geq\mu^+_{i(*)}$.
\item [$(e)$] For every $i<i(*)$ and every measurable $\mu_i>\lambda$, there exists a sequence of measurables $\langle\lambda_{i,n}:n<\omega\rangle$ so that ${\rm tcf}(\prod\limits_{n\in\omega}\lambda_{i,n},<_{J_\omega^{\rm bd}})=\mu_i$ and ${\rm tcf} (\prod\limits_{n\in\omega}\lambda_{i,n}^+,<_{J_\omega^{\rm bd}})=\mu_i^+$.
\end{enumerate}
\end{theorem}

\par\noindent\emph{Proof}.\newline 
Let $\mathbb{P}$ be the extender-based Prikry forcing with respect to $\lambda$ and $\mu_{i(*)}$. The forcing notion $\mathbb{P}$ satisfies the $\lambda^{++}$-cc, hence cardinals above $\lambda^+$ are preserved. It also satisfies the Prikry property, so by the completeness of the pure order $\leq^*$ we conclude that every cardinal below $\lambda$, as well as $\lambda$ itself, is preserved. It is also known that $\lambda^+$ is not collapsed, so $(a)$ holds.

Each Prikry sequence shows that $\cf(\lambda)=\aleph_0$ in the generic extension, hence $(b)$ holds. By the Prikry property (and the completeness of $\leq^*$) no bounded subsets of $\lambda$ are added, so if one begins with GCH in the ground model then $(c)$ holds in ${\rm \bf V}^{\mathbb{P}}$. We shall prove, herein, that the generic object adds a scale of $\mu_{i(*)}$-many Prikry sequences, hence we get $(d)$. 

Let $G\subseteq\mathbb{P}$ be a generic set and $\alpha<\mu_{i(*)}$. Set $G^\alpha=\bigcup\{p^\alpha:p\in G\}$. Given two ordinals $\alpha<\beta$ we shall prove that $G^\alpha(n)<G^\beta(n)$ for almost every $n\in\omega$.
Choose a condition $q\in G$ so that $\alpha,\beta\in{\rm supp}(q)$. Let $\gamma$ be ${\rm mc}(q)$, and set:

$$
A_\gamma=\{\nu<\kappa:\pi_{\gamma\alpha}(\nu)<\pi_{\gamma\beta}(\nu)\}.
$$

Since $\alpha<\beta$ we know that $A_\gamma\in U_\gamma$. Notice that $q^\alpha \subseteq G^\alpha$ and $q^\beta \subseteq G^\beta$, and without loss of generality $n_0=lg(q^\alpha)=lg(q^\beta)$. We shall prove that $n_0\leq n \Rightarrow G^\alpha(n)<G^\beta(n)$.

First we intersect $T^q_{q^{\rm mc}}$ at every level with $A_\gamma$. Using property $(4)$ of the order $\leq_{\mathbb{P}}$ we know that if $q\leq r$ then $r^\alpha\setminus q^\alpha=\pi_{\gamma\alpha}''(r^\gamma\setminus q^\gamma)$ and $r^\beta\setminus q^\beta=\pi_{\gamma\beta}''(r^\gamma\setminus q^\gamma)$.
Inasmuch as $r^\gamma\in T^q_{q^{\rm mc}}$ we have:

$$
\nu\in r^\gamma\setminus q^\gamma \Rightarrow \pi_{\gamma\alpha}(\nu)<\pi_{\gamma\beta}(\nu).
$$

which amounts to $G^\alpha(n)<G^\beta(n)$ whenever $n\geq n_0$. It follows that at least $\mu_{i(*)}^+$-many distinct Prikry sequences are added, so $2^\lambda\geq\lambda^{\aleph_0}\geq \mu_{i(*)}^+$. 

For proving $(e)$ fix any measurable $\mu_i$, suppose that it corresponds to the ordinal $\beta<\mu_{i(*)}^+$ and let $\langle\lambda_n:n\in\omega\rangle$ be the associated Prikry sequence through $U_\beta$. Assume $\name{t}$ is a name for a sequence in $\prod\limits_{n\in\omega}\lambda_n$. We have to show that there exists some ordinal $\alpha<\beta$ and a condition $q\in\mathbb{P}$ so that $q$ forces $t<_{J_{\rm bd}^\omega} \langle G^\alpha(n):n\in\omega\rangle$. As each Prikry sequence (in the diagonal Prikry forcing used above) adds a dominating family (see \cite{MR2768695}, Section 1.3) we are done.

\hfill\qedref{eeextender}

So far we know that each $\mu_i^+$ enters the spectrum of characters.
The last forcing that we shall employ is the L\'evy collapse from \cite{MR0268037}. Suppose $\lambda>\kappa=\cf(\kappa)$, and $\lambda$ is a regular cardinal which satisfies $\lambda^{<\kappa}=\lambda$. We define the forcing notion ${\rm Levy}(\kappa,\lambda)$. A condition in ${\rm Levy}(\kappa,\lambda)$ is a partial function $f:\kappa\rightarrow\lambda$ such that $|{\rm dom}f|<\kappa$. The order is inclusion.

${\rm Levy}(\kappa,\lambda)$ is $\lambda^+$-c.c. (provided that $\lambda^{<\kappa}=\lambda$) and $\kappa$-complete, hence cardinals below $\kappa$ or above $\lambda$ are preserved. The completeness of ${\rm Levy}(\kappa,\lambda)$ gives also preservation of true cofinalities, as demonstrated in the following:

\begin{claim}
\label{llevy} Preservation of tcf by the L\'evy collapse. \newline 
Suppose $\mathbb{P}$ is a $\lambda^+$-c.c. forcing notion, $\kappa=\cf(\lambda)<\lambda, \lambda<\chi<\mu$, and $\chi=\cf(\chi),\mu^{<\chi}=\mu$ in ${\rm \bf V}^{\mathbb{P}}$. \newline 
Let $\mathbb{R}$ be ${\rm Levy}(\chi,\mu)$. 
Assume further that $\bar{\lambda}$ is an increasing sequence of regular cardinals which tends to $\lambda$, and ${\rm \bf V}^{\mathbb{P}}\models {\rm tcf}(\prod\bar{\lambda},<_J)=\mu$ for some ideal $J$ on $\kappa$. 
Then ${\rm \bf V}^{\mathbb{P}\ast\mathbb{R}}\models {\rm tcf}(\prod\bar{\lambda},<_J)=\chi$. \newline 
Moreover, if $\mathbb{R}_i={\rm Levy}(\chi_i,\mu^+_i)$ for every $i<i(*)$, $\mathbb{R}=\prod\{\mathbb{R}_i:i<i(*)\}$ with Easton support, and ${\rm\bf V}^\mathbb{P} \models {\rm tcf}(\prod\bar{\lambda}_i,<_J)=\mu_i^+$ for every $i<i(*)$, then ${\rm \bf V}^{\mathbb{P}\ast\mathbb{R}}\models {\rm tcf}(\prod\bar{\lambda}_i,<_J) =\chi_i$.
\end{claim}

\par\noindent\emph{Proof}.\newline 
By the regularity of $\chi$ we know that $\mathbb{R}$ is $\chi$-complete. It follows that forcing with $\mathbb{R}$ does not add bounded subsets of $\chi$ to the universe. In particular, no new function in $(\prod\bar{\lambda}_i,<_J)$ is added, as $\lambda<\chi$.

since ${\rm\bf V}^\mathbb{P} \models {\rm tcf}(\prod\bar{\lambda}_i,<_J)=\mu$ and $\mu$ becomes $\chi$ in ${\rm \bf V}^{\mathbb{P}\ast\mathbb{R}}$ we have ${\rm \bf V}^{\mathbb{P}\ast\mathbb{R}}\models {\rm tcf}(\prod\bar{\lambda}_i,<_J) =\chi$. The same holds for the product $\mathbb{R}=\prod\{\mathbb{R}_i:i<i(*)\}$ with Easton support.

\hfill \qedref{llevy}

Recall that in Theorem \ref{oonestep}, $\kappa=\cf(\lambda)$ may be uncountable. In the main theorem below we insert many cardinals into ${\rm Sp}_\chi(\lambda)$, but the cofinality of $\lambda$ must be countable:

\begin{theorem}
\label{mt} The main theorem. \newline 
Suppose $\lambda$ is supercompact, and $\langle\mu_i:i<i(*)\rangle$ is a sequence of measurable cardinals above $\lambda$. Let $\langle\chi_i: i<i(*)\rangle$ be an increasing sequence of regular cardinals above $\lambda$ such that $\chi_i\leq\mu_i$ for every $i<i(*)$. \newline 
Then one can force $\{\chi_i: i<i(*)\} \subseteq {\rm Sp}_\chi(\lambda)$.
\end{theorem}

\par\noindent\emph{Proof}.\newline 
First observe that $\mu_i^+\in{\rm Sp}_\chi(\lambda)$ after forcing with the extender-based Prikry forcing $\mathbb{P}$ of \ref{eeextender}. Indeed, for each $\mu_i^+$ we have (in the forcing extension) a sequence $\bar{\lambda_i}$ so that ${\rm tcf}(\prod\limits_{n\in\omega}\lambda_n,<_{\rm J^{\rm bd}_\omega})=\mu_i$ and ${\rm tcf}(\prod\limits_{n\in \omega} \lambda_n^+,<_{\rm J^{\rm bd}_\omega})=\mu_i^+$. It follows from Claim \ref{oonestep} that $\mu_i^+$ is realized as a member of ${\rm Sp}_\chi(\lambda)$.

In purpose to incorporate the $\chi_i$-s we collapse the $\mu_i$-s. Without loss of generality, $\chi_i\leq\mu_i$ for every $i<i(*)$ (one has to choose for every $i<i(*)$ a measurable cardinal $\mu>\chi_i$ and rename the sequence of measurables).
For every $i<i(*)$ let $\mathbb{R}_i$ be ${\rm Levy}(\chi_i,\mu_i^+)$. Let $\mathbb{R}$ be the product $\prod\{\mathbb{R}_i:i<i(*)\}$ with Easton support. We claim that the following holds in ${\rm \bf V}^{\mathbb{P}\ast\mathbb{R}}$ for every $i<i(*)$:

$$
\exists\bar{\lambda}_i, {\rm tcf}(\prod\bar{\lambda}_i)\leq{\rm tcf}(\prod\bar{\lambda}_i^+)=\chi_i.
$$

This assertion follows from \ref{llevy}. Working in ${\rm \bf V}^{\mathbb{P}\ast\mathbb{R}}$, Claim \ref{oonestep} gives now $\chi_i\in{\rm Sp}_\chi(\lambda)$ for every $i<i(*)$, so we are done.

\hfill \qedref{mt}

The extender-based Prikry forcing can be used for violating SCH on $\aleph_\omega$. One has to interlace L\'evy collapses in the forcing conditions, as shown in \cite{MR2768695}. Moreover, the pcf structure (manifested in Theorem \ref{eeextender}) remains similar, in particular the Prikry sequences form a scale (see \cite{master}). However, the combinatorial argument of Claim \ref{oonestep} involves the normality of the ultrafilters, so we may ask:

\begin{question}
\label{qqquestion} Is it forcable that ${\rm Sp}_\chi(\aleph_\omega)\supseteq\{\chi_i:i<i(*)\}$ for every increasing sequence of regular cardinals $\langle\chi_i:i<i(*)\rangle$ above $\aleph_\omega$?
\end{question}

We conclude with another problem. As mentioned in the introduction, the above methods provide a tool for possessing a cardinal in ${\rm Sp}_\chi(\lambda)$. We are interested also in the other side of the coin:

\begin{question}
\label{qq} Assume $\lambda>\cf(\lambda)$. Is it possible that ${\rm Sp}_\chi(\lambda)$ is not a convex set?
\end{question}

\newpage 

\bibliographystyle{amsplain}
\bibliography{arlist}

\end{document}